\documentclass[12pt]{article}
\usepackage{amssymb,amsmath,amsthm}
\usepackage{epsfig}
\usepackage{graphicx}

\begin{document}
\centerline{Jozef H. Przytycki}
\centerline{\bf Progress in distributive homology: }
\centerline{from q-polynomial of rooted trees to Yang-Baxter homology}
\ \\ \ \\
{\bf Abstract:}\\
{\footnotesize This is an extended abstract of the talk given at the Oberwolfach Workshop
``Algebraic Structures in Low-Dimensional Topology", 25 May -- 31 May 2014.
My goal was to describe progress in distributive homology from 
the previous Oberwolfach Workshop June 3 - June 9, 2012, in particular my work on 
Yang-Baxter homology; however I concentrated my talk on my recent discovery 
of q-polynomial of a rooted tree; the appropriate topic as my talk was on May 30, 2014, the 
30 anniversary of the Jones polynomial, and the polynomial has its roots in the Kauffman bracket 
approach to the Jones polynomial.}
\ \\ \ \\ \ \\

We start with a long historical introduction beginning with Heinrich K\"uhn (1690-1769), 
Carl Leonhard Gottlieb Ehler (1685-1753), and 
Leonard Euler (1707-1783) and we argue 
that topology (geometria situs) started in Gda{\'n}sk (Danzig) about 1734. We mention the work of 
Celestyn Burstin (1888-1938) and Walter Mayer (1887--1948), (1929, distributivity) and Samuel Eilenberg (1913-1998) 
 (homological algebra). We complete the historical 
summary by celebrating 30 years of the Jones polynomial (May 30, 1984, V.F.R.Jones wrote a letter to 
J.Birman announcing his construction of a new link polynomial). Thus it is appropriate to describe 
today a new simple invariant of rooted trees.
Let $T$ be a plane rooted tree then $Q(T)\in Z[q]$ is defined by the initial condition $T(\bullet)=1$ and 
the recursion relation 
$$Q(T)= \sum_{v\in L(T)}q^{r(T,v)}Q(T-v), \mbox{ where $L(T)$ is the set of leaves of $T$,  }$$
and  $r(T,v)$ is the number of edges of $T$ to the right of the path connecting $v$ with the root $v_0$.
For example $Q(\bigvee)= (1+q)=[2]_q$ or more generally $Q(T_n)= [n]_q!$, where $T_n$ is a star with $n$ rays and 
$[n]_q=1+q+...+q^{n-1}$.\\ 
Theorem: Let $T_1 \vee T_2$ be the wedge (or root) product ({\psfig{figure=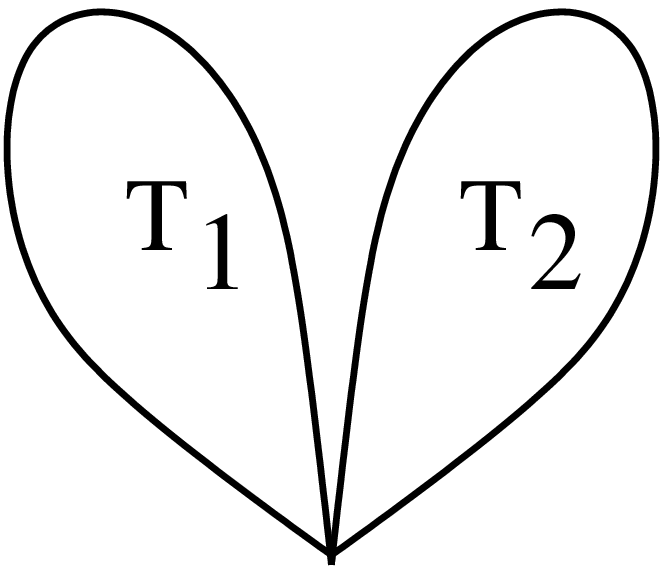,height=1.1cm}}). Then:
$$Q(T_1 \vee T_2)= \binom{E(T_1)+E(T_2)}{E(T_1)}_qQ(T_1)(Q(T_2)$$
Proof: We proceed by induction on $E(T)$, with obvious initial case of $E(T_1)=0$ or $E(T_2)=0$.
Let $T$ be a rooted plane tree with $E(T_1)E(T_2)> 0$, then we have:
$$Q(T)= \sum_{v\in L(T)}  q^{r(T,v)}Q(T-v)= $$
$$\sum_{v\in L(T_1)}  q^{r(T_1,v)+E(T_2) }Q((T_1-v)\vee T_2)+
\sum_{v\in L(T_2)}  q^{r(T_2,v)}Q(T_1\vee (T_2-v))  \stackrel{inductive\ assumption }{=}$$
$$\sum_{v\in L(T_1)}  q^{r(T_1,v)+E(T_2) }{E(T_1)+E(T_2)-1 \choose E(T_1)-1,E(T_2)}_q Q(T_1-v)Q(T_2) + $$
$$\sum_{v\in L(T_2)}  q^{r(T_2,v)}{E(T_1)+E(T_2)-1 \choose E(T_1),E(T_2)-1}_qQ(T_1)Q(T_2-v) =$$
$$ Q(T_2)q^{E(T_2)}{E(T_1)+E(T_2)-1\choose E(T_1)-1,E(T_2)}_q\sum_{v\in L(T_1)}  q^{r(T_1,v)}Q(T_1-v) +$$
$$ Q(T_1){E(T_1)+E(T_2)-1 \choose E(T_1),E(T_2)-1}_q\sum_{v\in L(T_2)}  q^{r(T_2,v)}Q(T_2-v) =$$
$$Q(T_1) Q(T_2)(q^{E(T_2)}{E(T_1)+E(T_2)-1 \choose E(T_1)-1,E(T_2)}_q +
{E(T_1)+E(T_2)-1\choose E(T_1),E(T_2)-1}_q)=$$
$$ Q(T_1) Q(T_2){E(T_1)+E(T_2)\choose E(T_1),E(T_2)}_q \mbox{ as needed }.$$

Corollary: 
\begin{enumerate}
\item[(i)] If a plane rooted tree is a wedge of $k$ trees
 (\psfig{figure=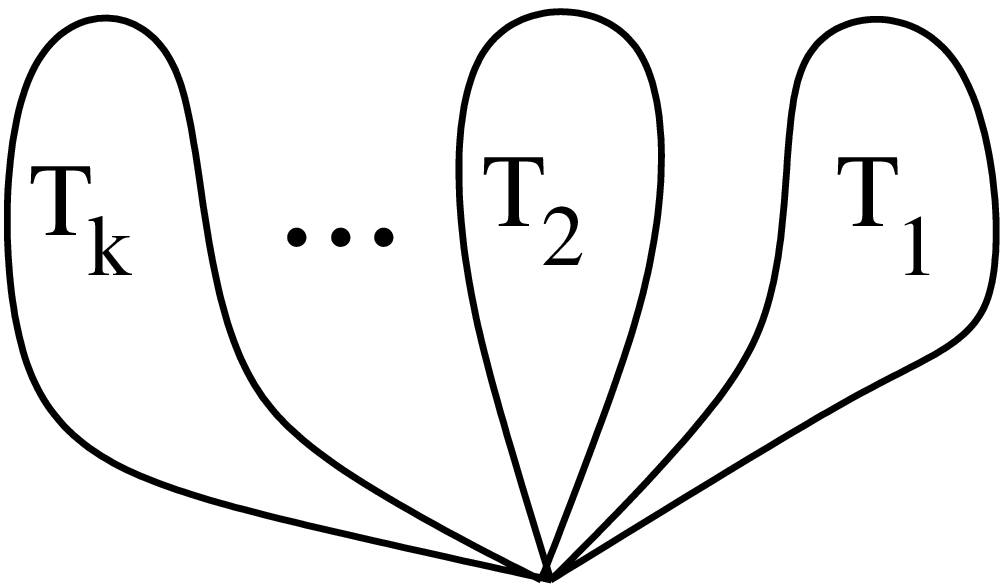,height=1.5cm}) and
$$T=T_{k} \vee ... \vee  T_2 \vee T_1,  \mbox{ then }$$
$$Q(T)= \binom{E_k+E_{k-1}+...+E_1}{E_k,E_{k-1},...,E_1}_qQ(T_k)Q(T_{k-1})\cdots Q(T_1).$$
where $E_i=|E(T_i)|$ is the number of edges in $T_i$.
\item[(ii)] (State product formula) 
$$ Q(T) = \prod_{v\in V(T)}W(v), $$
where $W(v)$ is a weight of a vertex (we can call it a Boltzmann weight) defined by:
$$W(v)= \binom{E(T^v)}{E(T^v_{k_v}),...,E(T^v_{1})}_q,$$
where $T^v$ is a subtree of $T$ with vertex $v$ (part of $T$ above $v$, in other words growing from $v$)  and $T^v$
can be decomposed into wedge of trees: $T^v= T^v_{k_v} \vee ... \vee  T^v_2 \vee T^v_{1}.$
\item[(iii)] (change of a root). Let $e$ be an edge of a tree $T$ with endpoints $v_1$ and
$v_2$ and $E_1$ be the number of edges on the $v_1$ part of the edge, and $E_2$ the number of edges of $T$ on
the $v_2$ side of $e$; 
$$\mbox{Thus $T=$\psfig{figure=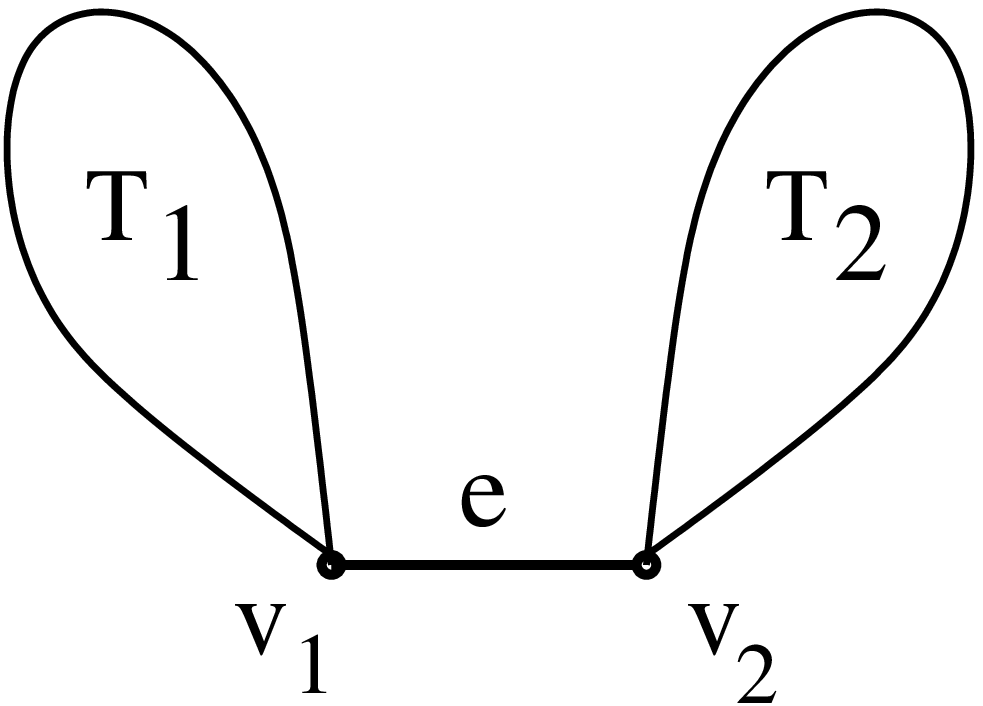,height=0.7cm} and then } 
Q(T,v_1)=\frac{[E_1+1]_q}{[E_2+1]_q}Q(T,v_2).$$
\end{enumerate}

\begin{proof}
(i) Formula z (i)  follows by using several times the formula
$$Q(T_2 \vee T_1)= {E(T_2)+E(T_1) \choose E(T_2),E(T_1)}_qQ(T_2)(Q(T_1),$$
as we have: 
$$\binom{a_k+a_{k-1}+...+a_2+a_1}{a_k,a_{k-1},...,a_2,a_1}_q=
\binom{a_{k-1}+...+a_2+a_1}{a_{k-1},...,a_2,a_1}_q \binom{a_k+a_{k-1}+...+a_2+a_1}{a_k,a_{k-1}+...+a_2+a_1}_q=...$$
$$= \binom{a_2+a_1}{a_2,a_1}_q\binom{a_3+a_2 + a_1}{a_3,a_2+a_1}_q\binom{a_4+a_3+a_2 + a_1}{a_4,a_3+a_2+a_1}_q...
\binom{a_k+a_{k-1}+...+a_2+a_1}{a_k,a_{k-1}+...+a_2+a_1}_q$$
(ii) Formula (ii) follows by using (i) several times.
\end{proof}

One can propose many modifications and generalizations of the polynomial $Q(T)$, for example, for a graph 
with a base point we can take the set (or the sum) over all spanning trees of $Q(T)$ but we propose below 
the one having close relation with knot theory.

Let $T$ be a plane rooted tree and $f:L(T) \to N$ a function from leaves of $T$ to positive integers. 
We call $f$ a delay function as our intuition is that a leaf with value $k$ cannot be used before $k$th move.
Formally $Q(T,f)$ is defined by recursive relation:
$$Q(T,f)= \sum_{v \in L_1(T)}q^{r(T,v)}Q(T-v,f_v),$$
where $ L_1(T)$ is a set of leaves for which $f$ is equal to $1$. $f_v(u)= max(1, f(u)-1)$ if $u$ is also 
a leaf of $T$, and it is equal to $1$ if it is a new leaf of $T-v$.\\
Example. For a rooted tree with delay function the polynomial $Q(T)$ is not necessary a product of 
cyclotomic polynomials, the simplest example is given by trees \psfig{figure=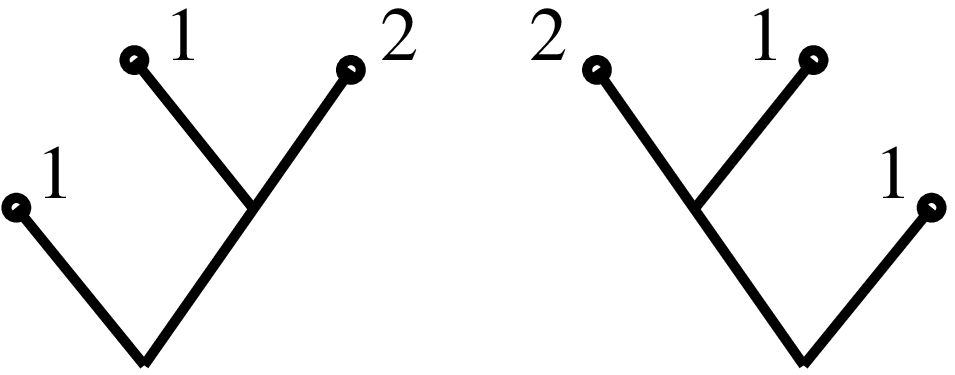,height=0.7cm}
 with polynomials equal respectively 
$q(1+q+2q^2+q^3) \mbox{ and } 1+2q + q^2+q^3.$ There is however one special situation when we can give a 
simple closed formula: Consider the ``delayed" tree $T=T_k^{s_k}\vee ...\vee T_2^{s_2} \vee T_1^{s_1}$ 
(\psfig{figure=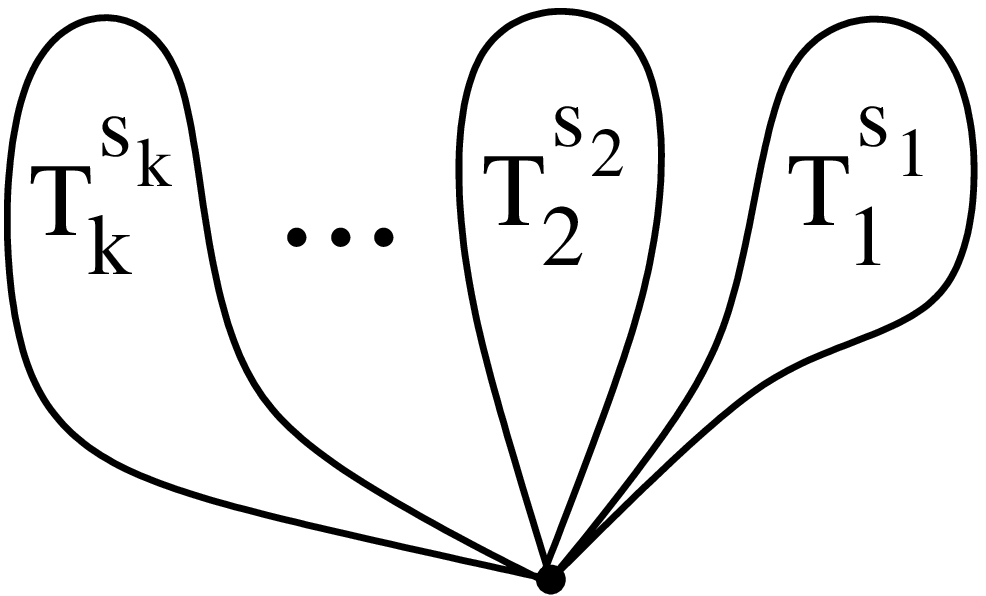,height=1.1cm}).
That is we assume that whole blocks have constant delay function (the block $T_i$ have leaves labelled $s_i$).
We assume also, for convenience, that $s_1=1$, $s_1 \leq s_2 \leq E_1 +1$, $s_2 \leq s_3 \leq E_2 + E_1 +1$,...,
$s_{k-1} \leq s_k \leq E_{k-1}+...+ E_2 + E_1 +1$ (here $E_i= |E(T_i)|$).
Then $$Q(T)=\binom{E_2+E_1 -s_2+1}{E_2,E_1-s_2+1}_q \binom{E_3+ E_2+E_1 -s_3+1}{E_3,E_2+ E_1-s_3+1}_q ...
\binom{E_k+...+E_1 -s_k+1}{E_k,E_{k-1}+...+E_1-s_k+1}_q$$
$$Q(T_1) Q(T_2)...Q(T_k).$$

We didn't reach yet relations neither  with knot theory nor with distributive structures; these should be left for the 
next occasion, however we finish the talk with one curious question and related observation.
Consider a chain complex over a commutative ring $k$ 
$$ {\mathcal C}:  ...\to C_{n+1} \stackrel{\partial_{n+1}}{\rightarrow}  C_n \stackrel{\partial_{n}}{\rightarrow}  C_{n-1}  
\stackrel{\partial_{n-1}}{\rightarrow}... \to  C_1 \stackrel{\partial_{1}}{\rightarrow} C_0 $$ 
and assume that $ {\mathcal C}$ comes from a presimplicial module $\partial_n=\sum_{i=0}^n(-1)^id_i$ 
where $0\leq i \leq n$,  
$d_id_j=d_{j-1}d_i$ for $i<j$. We ask whether it is useful (already used?) to consider q-version:
$C^q_n= C_n\otimes_k Z[q]$ and the q-map $\partial_n^q=\sum_{i=0}^nq^id_i$. Clearly $(C^q_n, \partial_n^q)$ is 
not generically a chain complex but we can make another use of it. For example, we can identify $x$ with $\partial^q(x)$,
that is to consider $ (\bigoplus_{n\geq 0}C^q_n)/(x-\partial^q(x))$.
Here an example which I learned from JP. Loday is very handy:\\
 Consider presimplicial set $(Y_n,d_i)$ where 
$Y_n$ is the set of topological rooted trees with $n$ ordered leaves (topological means that 
\psfig{figure=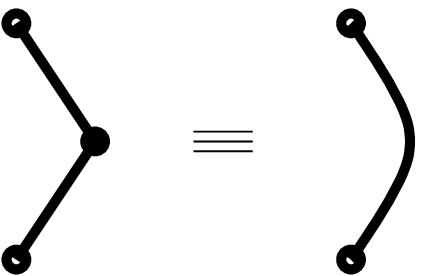,height=0.9cm}). We define $d_i(T)=T-v_i$, where $v_i$ is the $i$th leaf of $T$.
We can also introduce degeneracy maps $s_i:Y_i \to Y_{i+1}$ planting $\bigvee$ on the $i$th leaf. 
We check directly  that:\\
(1) $d_id_j=d_{j-1}d_i$ for $i<j$, \\
(2') $s_is_j = s_{j+1}s_i$ for $i<j$,\\
$ (3) \ \ \ d_is_j= \left\{ \begin{array}{rl}
 s_{j-1}d_i &\mbox{ if $i<j$} \\
s_{j}d_{i-1} &\mbox{ if $i>j+1$}
       \end{array} \right.
$\\
(
(4) $d_is_i=d_{i+1}s_i=Id $

The condition $s_is_i =s_{i+1}s_i$ does not hold (\psfig{figure=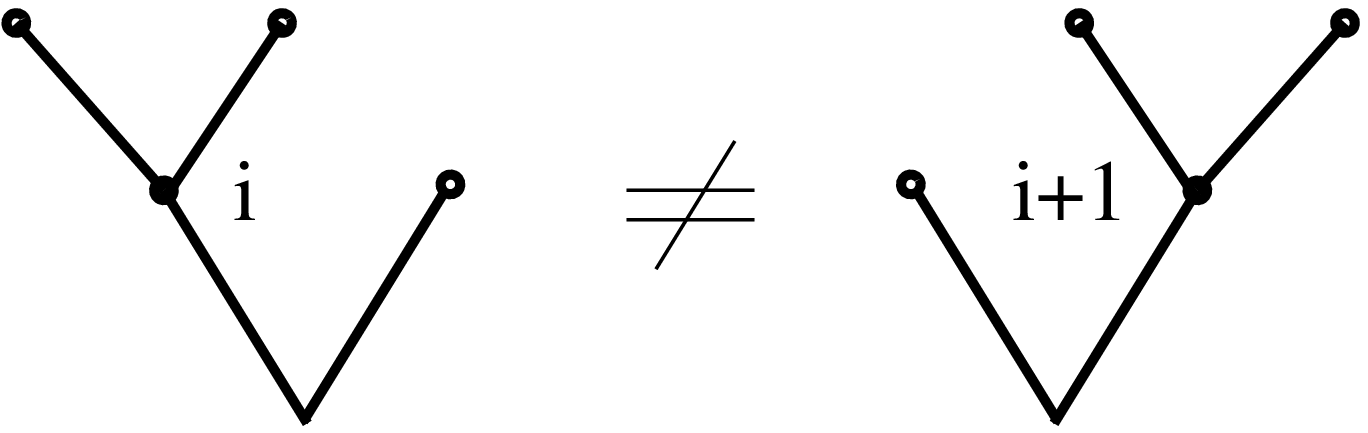,height=0.9cm}) 
so $(Y_n,d_i,s_i)$ is not a simplicial set but only an almost simplicial set.

Now consider the quotient of the sum $(\bigoplus_{n\geq 0}Z[q]Y_n)/(x-\partial^q(x))$.
It is a free $Z[q]$ module generated by $\bullet$ (tree without edges).We compute inductively that for 
a tree with $n$ leaves $T= [n]_q!\bullet$. It is not very sophisticated invariant so we can be glad that 
polynomial $Q(T)$ is more interesting.

Distributivity leads to another ``incomplete" simplicial set, this time condition (4) does not hold, 
but this should be put aside for the next report which will discuss also a generalization of 
distributive homology: -- Yang-Baxter homology. \\
\ \\
\centerline{\psfig{figure=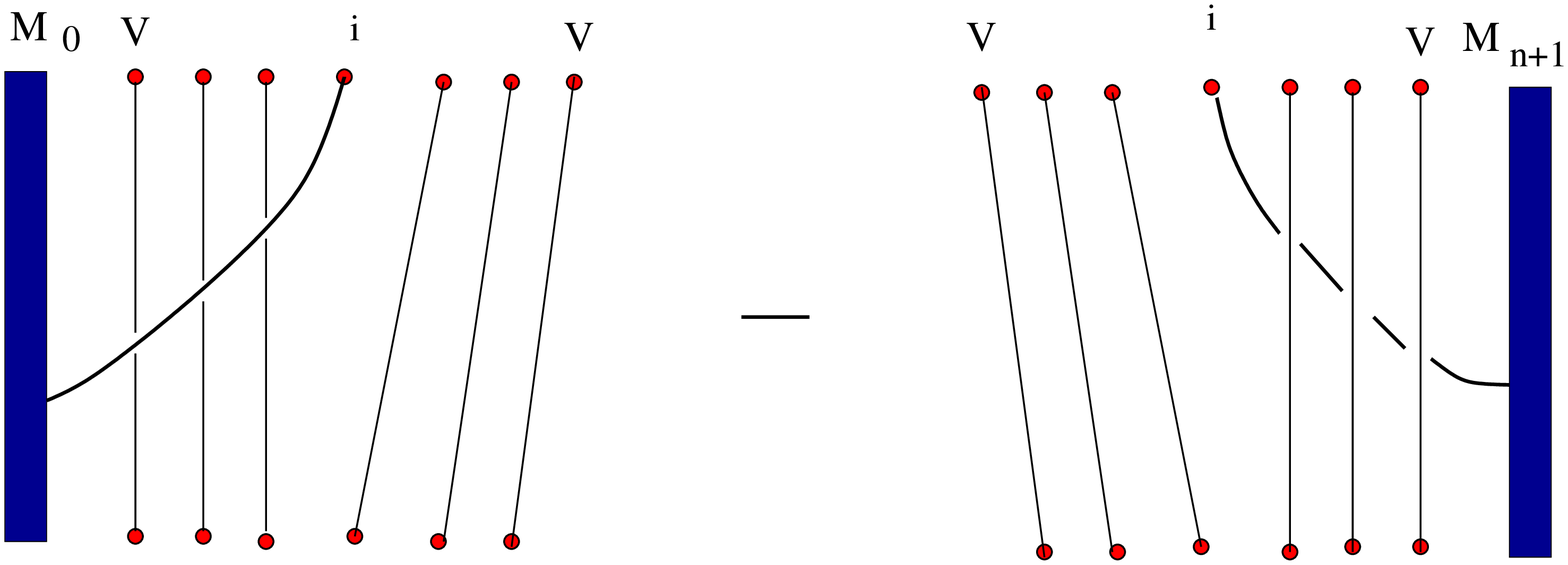,height=4.3cm}}
\ \\
\centerline{\bf \Large Thank You}

{\footnotesize George Washington University (przytyck@gwu.edu), UMD, and  University of Gda\'nsk}

\begin{thebibliography}{999}


\bibitem[CPP]{CPP} A.~Crans, J.~H.~Przytycki, K.~Putyra,
\textit{Torsion in one term distributive homology},
{\it Fundamenta Mathematicae}, \textbf{225}, May, 2014, 75-94.
e-print:\ {\tt  arXiv:1306.1506 [math.GT] }

\bibitem[DLP]{DLP}  M.~K.~Dabkowski, C.~Li, J.~H.~Przytycki, 
Catalan states of lattice crossing, preprint, February 2014.
(The polynomial $Q(T)$ was motivated by this paper and will be used in the follow up paper as an 
important  ingredient in analysis of the Kauffman bracket
of the lattice crossing \psfig{figure=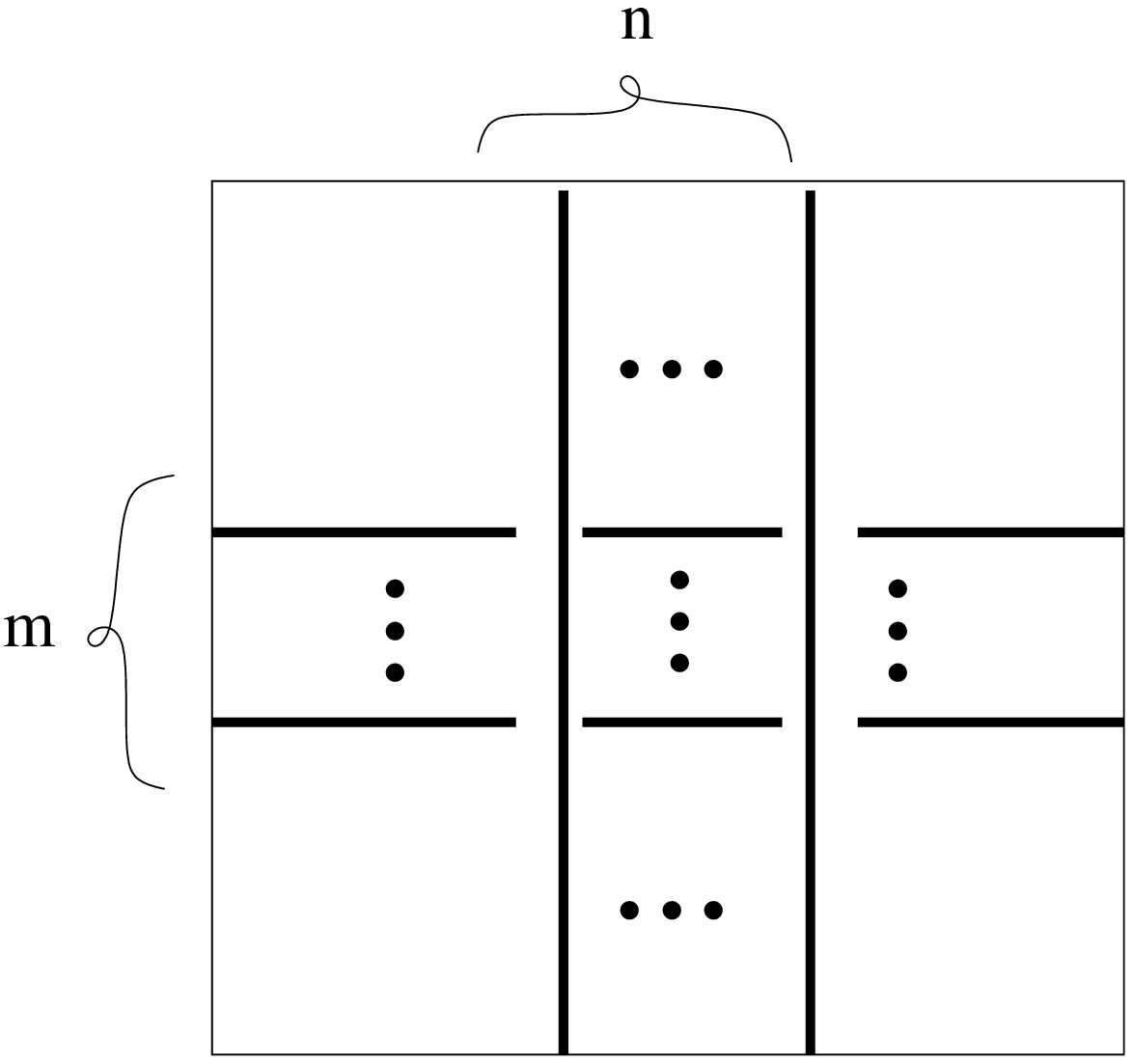,height=1.7cm}
.)

\bibitem [Prz-1]{Prz-1}
J.~H.~Przytycki,
Distributivity versus associativity in the homology theory of algebraic structures,
 {\it Demonstratio Math.}, \textbf{44}(4), December 2011, 823-869; \
 e-print:\ {\tt http://front.math.ucdavis.edu/1109.4850}

\bibitem[Prz-2]{Prz-2}
J.~H.~Przytycki,
{\it Knot Theory and related with knots distributive structures;
Thirteen Gdansk Lectures}, Gdansk University Press, in Polish, June, 2012, pp. 115 (second, extended, edition,
in preparation).

\bibitem[Prz-3]{Prz-3}
J.~H.~Przytycki,
Knots and distributive homology, Chapter in: {\it New Ideas in Low Dimensional Topology}, to appear
(Editors: L.H.Kauffman, V.Manturov).

\bibitem[P-P]{P-P} J.~H.~Przytycki, K.~Putyra,
\textit{Homology of distributive lattices},
the {\it Journal of homotopy and related structures}, Volume \textbf{8}(1), 2013, pages 35-65;
e-print: \ {\tt arXiv:1111.4772 [math.GT]}

\bibitem[P-R]{P-R}
J.~H.~Przytycki,, W.~Rosicki, \textit{Cocycle invariants of
codimension 2 embeddings of manifolds},  {\it Banach Center Publications};
Recommended for publication, January, 2014; to appear December 2014.
e-print: \ {\tt   arXiv:1310.3030 [math.GT]}

\bibitem[P-S]{P-S} J.~H.~Przytycki, A.~S.~Sikora,
\textit{Distributive Products and Their Homology},
 {\it Communications in Algebra}, \textbf{42}(3), 2014, 1258-1269;
e-print: \ {\tt arXiv:1105.3700 [math.GT]}\\

\end{thebibliography}
\end{document}